\documentclass[11pt]{article}
\usepackage{amsmath, amsfonts, amsthm, latexsym, hyperref }
\usepackage{graphicx,pxfonts,pstricks}

\newtheorem{thm}{Theorem}

\newcommand{\N}{\mathbb{N}}
\newcommand{\Z}{\mathbb{Z}}

\newcommand{\Prob}{\mathbb{P}}

\makeatletter
\newcommand{\eqnum}{\leavevmode\hfill\refstepcounter{equation}\textup{\tagform@{\theequation}}}
\makeatother

\title{Monotonicity of Avoidance Coupling on $K_N$}
\author{Ohad N. Feldheim\thanks{Institute of Mathematics and its Applications,
University of Minnesota,
College of Science and Engineering,
Minneapolis, MN U.S.A. 55455 Email:ohad\_f@netvision.net.il.
This research was supported in part by the Institute for Mathematics
and its Applications with funds provided by the National Science
Foundation.}}

\begin{document}
\maketitle
\begin{abstract}
Answering a question by Angel, Holroyd, Martin, Wilson and Winkler in
\cite{AHMWW}, we show that the maximal number of non-colliding coupled simple random
walks on the complete graph $K_N$, which take turns, moving one at a time, is
monotone in $N$. We use this fact to couple $\lceil \frac N4 \rceil$ such walks on $K_N$, improving the previous
$\Omega(N/\log N)$ lower bound of Angel et al. We also introduce a new generalization of simple avoidance coupling which we call partially ordered simple avoidance coupling
and provide a monotonicity result for this extension as well.
\end{abstract}

\section{Introduction}
Let $G=([N],E)$ be a graph whose vertices are the set of integers $[N]=\{1,\dots,N\}$.
A \emph{simple random walk} on this graph is a Markov chain $(X_t)_{t\in\Z}$ of
elements in $[N]$ such that for all $t\in\Z$ the distribution of $X_t$ is
uniform on the neighbors of $X_{t-1}$.

A \emph{Simple Avoidance Coupling (SAC)} of $k$ walks on $G$ is a sequence of
random maps $(U_t)_{t\in \Z}$ from $[k]$ to $[N]$ which
satisfy two conditions:
\begin{equation}\label{eq: SAC 1st}
\forall i\in [k] : (U_t(i))_{t\in\Z}\text{ is a simple random walk on }G
\end{equation}
\begin{equation}\label{eq: SAC 2nd}
\forall t\in \Z,\ 1 \le i<j\le k\ :\  \Prob\Big(U_t(i)=U_t(j)\Big)= \Prob\Big(U_t(i)=U_{t-1}(j)\Big)=
   0
\end{equation}

Angel, Holroyd, Martin, Wilson and Winkler introduce this notion in
\cite{AHMWW} in order to investigate couplings of $k$ simple random walks which move in turns
in discrete time and avoid collision.

One possible application of SACs on the complete graph $K_N$ is
semi-synchronous orthogonal frequency hopping. A communication network consists of several transmitters.
As there are overlaps between the transmission ranges they wish to use distinct frequencies at every given time. Malicious adversaries, each located in the vicinity of one of these transmitters, are trying to interfere with the communication by noising several frequencies at every given time.
Once an adversary hits his target transmitter's frequency he can tell that his interruption succeeded.
In order to avoid persistent interference the transmitters wish to change frequencies often.
Being unable to perfectly synchronize their clocks, the transmitter must take turns at hopping.
In this scenario it is desirable for each transmitter to perform a simple random walk as this would make each of its frequency changes (hops) independent from the past with maximal entropy. Independence is desirable since
 the adversary has some access to the frequency history of its target transmitter. An ideal hopping scheme in this setting is a SAC.

An important result of \cite{AHMWW} is that there exists a SAC of
$\Omega(N/\log N)$ walks on $K_N$. The authors also show in \cite[Theorem 6.1]{AHMWW} that when
$N=2^\ell+1$ for some $\ell\in \N$, there exists an avoidance coupling of
$2^{\ell-1}$ walks on $K_N$. Angel, Holroyd, Martin, Wilson and Winkler ask:
does the existence of an avoidance coupling of $k$ walks on $K_{N-1}$, imply the
existence of an avoidance coupling of $k$ walks on $K_{N}$. Our main result is a positive answer to this question:

\begin{thm}\label{thm: main1}
If there exists a simple avoidance coupling of $k$ walks on $K_{N-1}$, then there
exists a simple avoidance coupling of $k$ walks on $K_{N}$.
\end{thm}

Combining this with \cite[Theorem 6.1]{AHMWW}  we draw the following improved bound.

\begin{thm}\label{thm: main3}
There exists a simple avoidance coupling of $\lceil N/4 \rceil$ walks on $K_{N}$.
\end{thm}

We find it interesting that a byproduct of the proof of Theorem~\ref{thm: main1} is that in the extended coupling
on $K_{N+1}$ one find another $k+1$-th special walk, which is a simple random walk as well, but does not obey the order in which the walkers move in every round. In Section~\ref{sec: posac} we investigate this observation and discuss a possible extensions of avoidance coupling to models where the order by which the walks move changes from one round to the next, subject to some restrictions.

\subsection{Markovian Couplings}

In \cite{AHMWW} the authors give special attention
to Markovian Simple Avoidance Couplings. These have the property that whenever a walker's turn to move arrives,
he needs only to look at the current configuration walkers to determine the distribution of its next location.
In particular the simple avoidance coupling of $\Omega(N \log N)$ walkers on $K_N$ constructed in \cite{AHMWW} has this property, as does the coupling of $2^{\ell-1}$ walkers on $K_{2^\ell+1}$. While our extension
theorem does not preserve this property, we preserve the following weaker version. Consider a SAC in which each site of the underlying graph $K_{N}$ is assigned a label. At the end of every round
a random permutation is applied to these labels. Such a SAC is called \emph{Label Markovian} if whenever a walker's turn to move arrives he needs only to look at the current configuration of the walkers and, in addition, at the current labels of the vertices. Observe that every Markovian SAC
is also a Label Markovian SAC. It is straightforward to check that our
construction preserves Label Markov property.

\section{Background}
Probabilistic coupling of several stochastic processes sharing the same distribution, has been introduced to probability theory mainly as a tool to study and prove various properties of that common distribution. Such methods have been successfully used in showing properties such as monotonicity, stochastic dominance and convergence.

Nevertheless, probabilistic coupling can also be a subject of study. In this context, the natural question is "in what sense a collection of coupled identically distributed stochastic processes, is different from a collection of independent processes with the same distribution?". A classical example is that of two random walks on some finite graph $G$. If two independent random walks move on $G$, then they are bound to collide with high probability after a polynomial number of steps. Collisions occur even if a scheduler is allowed to control the times in which each walk makes his move (see \cite{CTW},\cite{TW}), and can be avoided only if the scheduler has some knowledge of the future of each walk, and only on special graphs (see \cite{G1}). On the other hand, there exist many graphs on which coupled random walks can easily avoid each other. On the cycle graph $C_n$ for example, two walks which start on non-adjecent vertices can avoid each other by moving in the same direction at every step - either clockwise or counter-clockwise. Coupling of walks on $K_N$, the complete graph on $N$ vertices, appears to be more difficult. In \cite{AHMWW}, the authors use various techniques inspired by discrete harmonic analysis to create an avoidance coupling of $\Omega(N/\log N)$ walks on $K_N$ and of $N/2-1$ walks for an infinity collection of special $N$-s. They also investigate avoidance coupling on $K^*_N$, the complete graph with loops on $N$ vertices, and obtain a lower bound of $N/4$ walks on this graph. The authors further show that no coupling exists for $N-1$ walks on $K^*_N$, if $N\ge 4$.

The research of avoidance couplings is closely related to that of Brownian motions which keep at least constant distance from each other. This subject and its relation to pursuit-evasion problems is investigated in \cite{BBC}, \cite{BBK} and \cite{K1}.

\section{Extending an avoidance coupling}

This section consists of the proof of Theorem~\ref{thm: main1}. Let $\mathcal{U}^{N-1}_k = \big(U_t(j)\big)_{t\in\Z,j\in[k]}$ be a SAC of $k$ walks on $K_{N-1}$. Our goal is to define $\mathcal{W}^{N}_k$, a SAC of $k$ walks on $K_{N}$.

\subsection{The extended coupling}\label{sec: ext}

We begin by introducing an auxiliary sequence of random permutations.
Let $P_0\in S_{N}$ be a uniformly chosen random permutation in $S_{N}$. Let $(a_t)_{t\in\Z}$ be an i.i.d. sequence where
$a_0$ is a uniformly chosen element of $[N-1]$. For $t\in \N$ define inductively $P_t,P_{-t}\in{S_N}$ as follows.
\begin{align*}
P_t&:=P_{t-1}\circ (N\ a_t ),\\
P_{-t}&:= P_{1-t} \circ (N\ a_{1-t}),
\end{align*}
where $(a\ b)$ is the transposition of the two elements $a$ and $b$.

Write $\mathcal{P}^{N}=(P_t)_{t\in\Z}$. It is straightforward to check that $\mathcal{P^{N}}$ is a stationary Markov chain on $S_{N}$ which is independent from $\mathcal{U}^{N-1}_k$.

We define
$\mathcal{W}^{N}_k=\big(W_t(j)\big)_{t\in\Z,j\in[k]}$ where $W_t:[k]\to[N]$, as follows:
\begin{equation*}
W_t(j) = P_tU_t(j),    \quad j\in [k], \: t\in \Z.
\end{equation*}

An example of $\mathcal{U}^{5}_2$, $\mathcal{P}^{6}$ and $\mathcal{W}^{6}_2$ is given in Figure~1.
Below we prove that $\mathcal{W}^N_k$ is an avoidance coupling of $k$ walks on $K_{N}$.

%In this section we define a \emph{vertex extension} (VE) of a
%SAC of $k$ walks on $K_{N}$. We then prove that the first coordinate of this VE is a
%SAC of $k$ walks on $K_{N+1}$, concluding the proof of Theorem~\ref{thm:
%main2}.

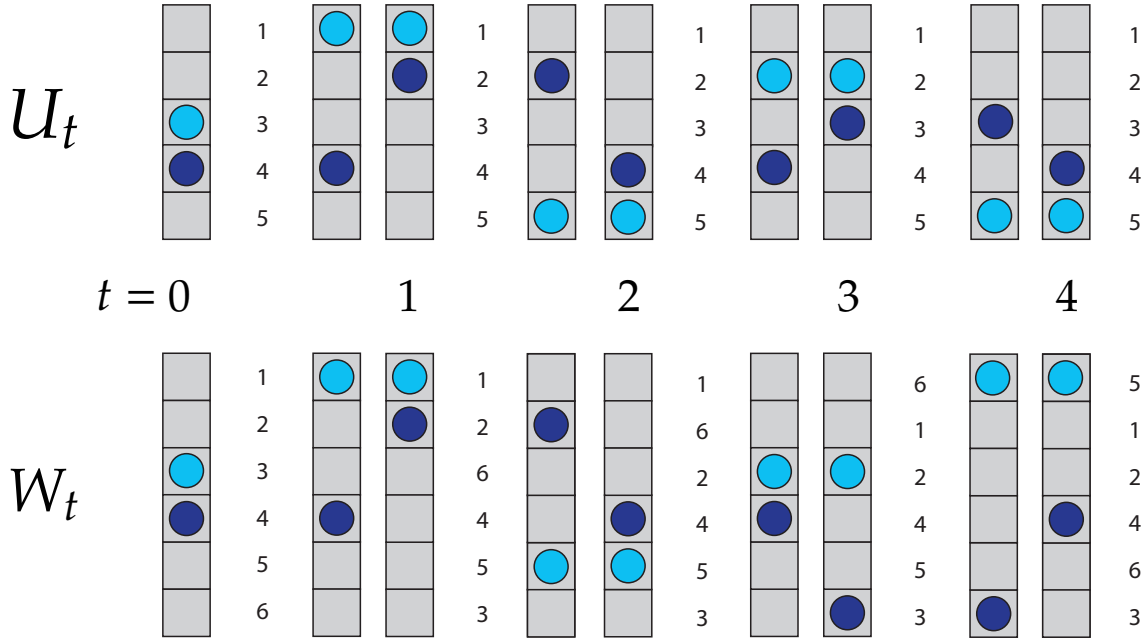
\begin{figure}[ht!]
  \centering
    \includegraphics[width=\textwidth]{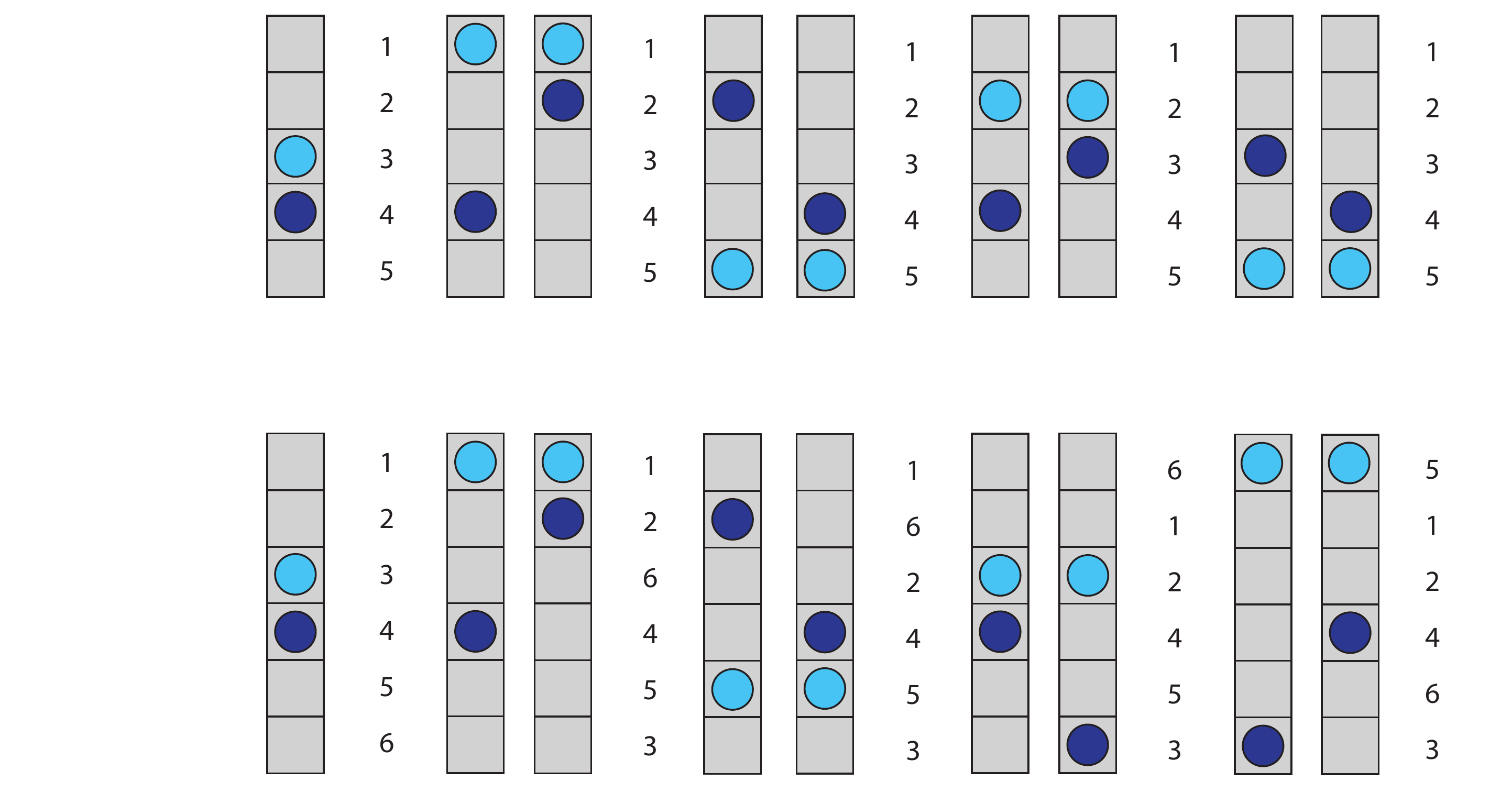}
  \rput(-7,7.8){\resizebox{1cm}{!}{$U_t$}}
  \rput(-7,2.8){\resizebox{1cm}{!}{$W_t$}}
  \rput(-5.9,5.45){\resizebox{0.8cm}{!}{$t=$}}
  \rput(-5.2,5.45){\resizebox{0.3cm}{!}{$0$}}
  \rput(-2.15,5.45){\resizebox{0.3cm}{!}{$1$}}
  \rput(0.78,5.45){\resizebox{0.3cm}{!}{$2$}}
  \rput(3.7,5.45){\resizebox{0.3cm}{!}{$3$}}
  \rput(6.6,5.45){\resizebox{0.3cm}{!}{$4$}}
  \caption{Above: $U_t$, a SAC of $2$ walks on $K_5$. Below: $W_t$, the extended SAC on $K_6$. The label permutation $P_t$ is given at the end of every time unit. Observe that the light blue walk always moves before the dark one. Also observe how $W_t$ is determined by  $P_t$ and $U_t$.}
\end{figure}

\subsection{The extension is a SAC}

%Our purpose is to find $\sigma_t$ satisfying \eqref{eq: PSAC 1st}, and verify
%that $W$ satisfies \eqref{eq: PSAC 2nd} and \eqref{eq: PSAC 3rd} with respect
%to $\sigma_t$.
To show that $\mathcal{W}^{N}_k$ is a SAC we must show that is satisfies \eqref{eq: SAC 1st} and \eqref{eq: SAC 2nd}.
We begin by showing $\eqref{eq: SAC 2nd}$.

Let $t\in \Z$,$\ 1 \le i<j\le k$. We have

$$\Prob\Big(W_t(i)=W_t(j)\Big)=
\Prob\Big(P_tU_t(i)=P_tU_t(j)\Big)=
\Prob\Big(U_t(i)=U_t(j)\Big)=0$$
Where the central equality uses the fact that $P_t$ is a permutation and the
right-most equality follow from the fact that $\mathcal{U}^{N-1}_k$ satisfies $\eqref{eq: SAC 2nd}$.

Recall the definition of the sequence $(a_t)_{t\in\Z}$
and write $P'_t$ for the transposition $(N\ a_t)$. We have
\begin{equation}
\begin{aligned}
\Prob\Big(W_t(i)=W_{t-1}(j)\Big)&=
\Prob\Big(P_{t}U_{t}(i)=P_{t-1}U_{t-1}(j)\Big)=
\Prob\Big(P_{t-1}\circ P'_t)U_{t}(i)=P_{t-1}U_{t-1}(j)\Big)\\
&=\Prob\Big(P'_tU_{t}(i)=U_{t-1}(j)\Big)=0,
\end{aligned}\label{eq: SAC 2n pt 2}
\end{equation}
where the the last equality follows from the fact that $\mathcal{U}^{N-1}_k$ satisfies $\eqref{eq: SAC 2nd}$,
and from the fact that $U_t(i),U_{t-1}(j)\in [N-1]$.

We are left with showing that $\mathcal{W}^{N}_k$ satisfies \eqref{eq: SAC 1st}.
Fix $j\in [k]$, we must show that $W_t(j)$ is a simple random walk on $K_{N}$.
Equivalently -- for every $\ell\in\N$, every
history $w_{t-\ell},...,w_{t-1}\in [N]$ such that
$$\Prob\Big(W_{t-1}(j)=w_{t-1},\dots,W_{t-\ell}(j)=w_{t-\ell}\big)>0,$$
and for every $v \neq w_{t-1}$, we have
\begin{equation}\label{eq: showing simple random walkness1}
\Prob\Big(W_t(j)=v\ \big|\ W_{t-1}(j)=w_{t-1},\dots,W_{t-\ell}(j)=w_{t-\ell}\big)=\frac1{N-1}.
\end{equation}

To obtain this we show a stronger claim. Fix $\ell\in\ N$ and let $p = (p_{t-\ell},...,p_{t-1})\in
(S_{N})^\ell$, $u= (u_{t-\ell},...,u_{t-1})\in [N+1]^\ell$.
Consider the event
$$A_t^{p,u}=\Big\{U_{t-1}(j)=u_{t-1},\dots,U_{t-\ell}(j)=u_{t-\ell}\text{ and
}P_{t-1}(j)=p_{t-1},\dots,P_{t-\ell}(j)=p_{t-\ell}\Big\}.$$

We show that for all $p,u$ such that $\Prob(A_t^{p,u})\neq0$ and for all
$v\neq p_{t-1}(u_{t-1})$ we have
\begin{equation}\label{eq: showing simple random walkness2}
\Prob\Big(W_t(j)=v\ \big|\ A_t^{p,u}\big)=1/N.
\end{equation}

Indeed, \eqref{eq: showing simple random walkness2} is stronger than
\eqref{eq: showing simple random walkness1}, as the values of
$P_{t-1},\dots,P_{t-\ell}$ and $U_{t-1}(j),\dots,U_{t-\ell}(j)$ determine the
values of $W_{t-1}(j),\dots,W_{t-\ell}(j)$.

Since $w_{t-1}=p_{t-1}(u_{t-1})\neq p_{t-1}(N)$ and using the fact that by \eqref{eq: SAC 2nd} we have
$$\sum_{n\in [N]\setminus w_{t-1}} \Prob\Big(W_t(j)=n\ \big|\ A_t^{p,u}\Big)=1,$$
it would suffice to show \eqref{eq: showing simple random walkness2} in the
case $v\neq p_{t-1}(N)$.
Thus, let $v\in [N]\setminus \{p_{t-1}(N),w_{t-1}\}$ and
use the total probability formula to write
\begin{align}\label{eq: main}
  \Prob\Big(W_t(j)=v\ \big|\ A_t^{p,u}\Big)
    &=\Prob\Big(W_t(j)=v \ \big|\ A_t^{p,u}, P_t(N)=v\Big)\Prob\Big(P_{t}(N)=v\Big) \notag \\
  &\ \ \ \ +\ \Prob\Big(W_t(j)= v\ \big|\ A_t^{p,u},  P_t(N)\neq v\Big)\Prob\Big(P_t(N)\neq v\Big) \notag \\
  &=\Prob\Big(U_t(j)=N \ \big|\ A_t^{p,u}, P_t(N)=v\Big)\cdot \frac 1 {N-1} \notag \\
  &\ \ \ \ +\Prob\Big(U_t(j)=P_t^{-1}(v) \ \big|\ A_t^{p,u}, P_t(N)\neq v\Big)\cdot \frac {N-2}{N-1}\notag\\
  &=0\cdot \frac 1 {N-1} +\Prob\Big(U_t(j)=P_t^{-1}(v) \ \big|\ A_t^{p,u}, P_t(N)\neq v\Big)\cdot \frac {N-2}{N-1}.
\end{align}
We now observe that
\begin{align}\label{eq: conditioned 1 over N-1}
&\Prob\Big(U_t(j)=P_t^{-1}(v) \ \big|\ A_t^{p,u}, v \neq P_{t}(N)\}\Big) = \notag \\
&\Prob\Big(U_t(j)=p_{t-1}^{-1}(v) \ \big|\ A_t^{p,u}, v \neq P_{t}(N)\}\Big) = \frac{1}{N-2}.
\end{align}
%\begin{align}\label{eq: conditioned 1 over N-1}
%&\Prob\Big(U_t(j)=P_t^{-1}(v) \ |\ A, v \notin\{P_{t}(N+1), P_{t-1}(N+1),p_{t-1}(u_{t-1})\}\Big) = \notag \\
%&\Prob\Big(U_t(j)=p_{t-1}^{-1}(v) \ |\ A, v \notin\{P_{t}(N+1), P_{t-1}(N+1),p_{t-1}(u_{t-1})\}\Big) = \notag \\
%&\Prob\Big(U_t(j)=p_{t-1}^{-1}(v) \ |\ v \notin\{P_{t}(N+1), p_{t-1}(u_{t-1})\}\Big) = \notag \\
%&\Prob\Big(U_t(j)=p_{t-1}^{-1}(v) \ |\ p_{t-1}^{-1}(v)\neq U_{t-1}(j)\Big) = \frac{1}{N-1}.
%\end{align}
where the first equality follows from the fact that for all $v
\notin\{P_{t}(N), P_{t-1}(N)\}$, we have $P_t^{-1}(v)=P_{t-1}^{-1}(v)$, and the last equality uses our assumption that $v\neq w_{t-1} = P_{t-1}U_{t-1}(j)$.

Plugging \eqref{eq: conditioned 1 over N-1} into \eqref{eq: main} we deduce
\eqref{eq: showing simple random walkness2}, concluding the proof. \qed

\section{Partially ordered avoidance coupling}\label{sec: posac}

Consider the following generalization of an avoidance coupling.
Let $R$ be a partial order on $[k]$.
An $R$ \emph{Partially Ordered Avoidance Coupling} (POSAC) of $k$ walks on $G$ is a sequence of random maps
\begin{equation*}
U_t:[m] \to [N], t\in \Z,
\end{equation*}
such that there exists a sequence of permutations $\sigma_t\in S_m$ which respect $R$ (i.e., $i<_Rj \rightarrow \sigma(i)<\sigma(j)$)
such that $(U_t)_{t\in \Z}$  and $\sigma_t$ satisfy two conditions:

\begin{enumerate}
\item  $\forall i\in [m] : (U_t[i])_{t\in\Z}\text{ is a simple random walk on }G,$ \eqnum\label{eq: PSAC 1st}
\item  $\forall t\in \Z,\ 1 \le \sigma_t(i)<\sigma_t(j)\le m\ :\ \Prob\Big(U_t(i)=U_{t-1}(j)\Big)=
    \Prob\Big(U_t(i)=U_t(j)\Big)=0$. \eqnum\label{eq: PSAC 2nd}
\end{enumerate}

A POSAC is a generalization of a SAC to a situation where the order in which the walks take turns can change from one round to the next, restricted by some partial order constraint (in the application to orthogonal hoping consider a situation where two transmitters can alter the order of their hops only if they receive each other's signal).

The proof of Theorem~\ref{thm: main1} extends in this case to the following.
\begin{thm}\label{thm: main2}
If there exists an $R$ POSAC of $k$ walks on $K_{N-1}$, then there
exists an $R$ POSAC of $k+1$ walks on $K_{N}$.
\end{thm}

Observe that in this case, although the extension does not allow adding additional relations it does allow increasing the number of walks.

\subsection{Extending a POSAC}

This section is dedicated to the proof of Theorem~\ref{thm: main2}.
Let $R$ be a partial order on $[k]$, let $\mathcal{U}^{N-1}_{k,R}$ be an $R$ POSAC of $k$ walks on $K_{N-1}$,
and let $(s_t)_{t\in\Z}$ be a sequence of permutations which respect $R$ and satisfy \eqref{eq: PSAC 2nd}.

Let $\mathcal{P}^{N}=(P_t)_{t\in\Z}$ and $\mathcal{W}^{N}_{k}=\big(W_t(j)\big)_{t\in\Z,j\in[k]}$ be as in section \ref{sec: ext}, and define
$\mathcal{W}^{N}_{k+1,R}=\big(W_t(j)\big)_{t\in\Z,j\in[k+1]}$ with $W_t(k+1):=P_t(N)$.

Observe that given $t\in \Z$, for any distinct $i,j\in [k+1]$ we have
\begin{equation}\label{eq: first part uniqueness ext}
\Prob\Big(W_t(i)=W_t(j)\Big)=
\Prob\Big(P_tU_t(i)=P_tU_t(j)\Big)=
\Prob\Big(U_t(i)=U_t(j)\Big)=0,
\end{equation}
as before. Using this we define $(\sigma_t)_{t\in\Z}$ in the following way. If there exists $b\in[k]$ such
that $W_{t-1}(b)=W_t(k+1)$ we set
\begin{equation}\label{eq: sigma case1}
\sigma_t (j)=\begin{cases} s_t(j)   & s_t(j)\le s_t(b) \\
                            s_t(b)+1 & j=m+1 \\
                            s_t(j)+1 & s_t(j) > s_t(b)
\end{cases}
\end{equation}
while otherwise we set
\begin{equation}\label{eq: sigma case2}
\sigma_t (j)=\begin{cases} s_t(j)+1 & j\le m \\
                             1    & j=m+1\ \ \ .
\end{cases}
\end{equation}

Our purpose is to show that $\mathcal{W}^{N}_k$ and $(\sigma_t)_{t\in\Z}$ satisfy \eqref{eq: PSAC 1st} and \eqref{eq: PSAC 2nd}. An example of $\mathcal{U}^{5}_{2,R}$, $\mathcal{P}^{6}$, $\mathcal{W}^{6}_{3,R}$ and $(\sigma_t)_{t\in\Z}$ is given in Figure~2.

%In this section we define a \emph{vertex extension} (VE) of a
%SAC of $k$ walks on $K_{N}$. We then prove that the first coordinate of this VE is a
%SAC of $k$ walks on $K_{N+1}$, concluding the proof of Theorem~\ref{thm:
%main2}.
\begin{figure}[ht!]
  \centering
    \includegraphics[width=\textwidth]{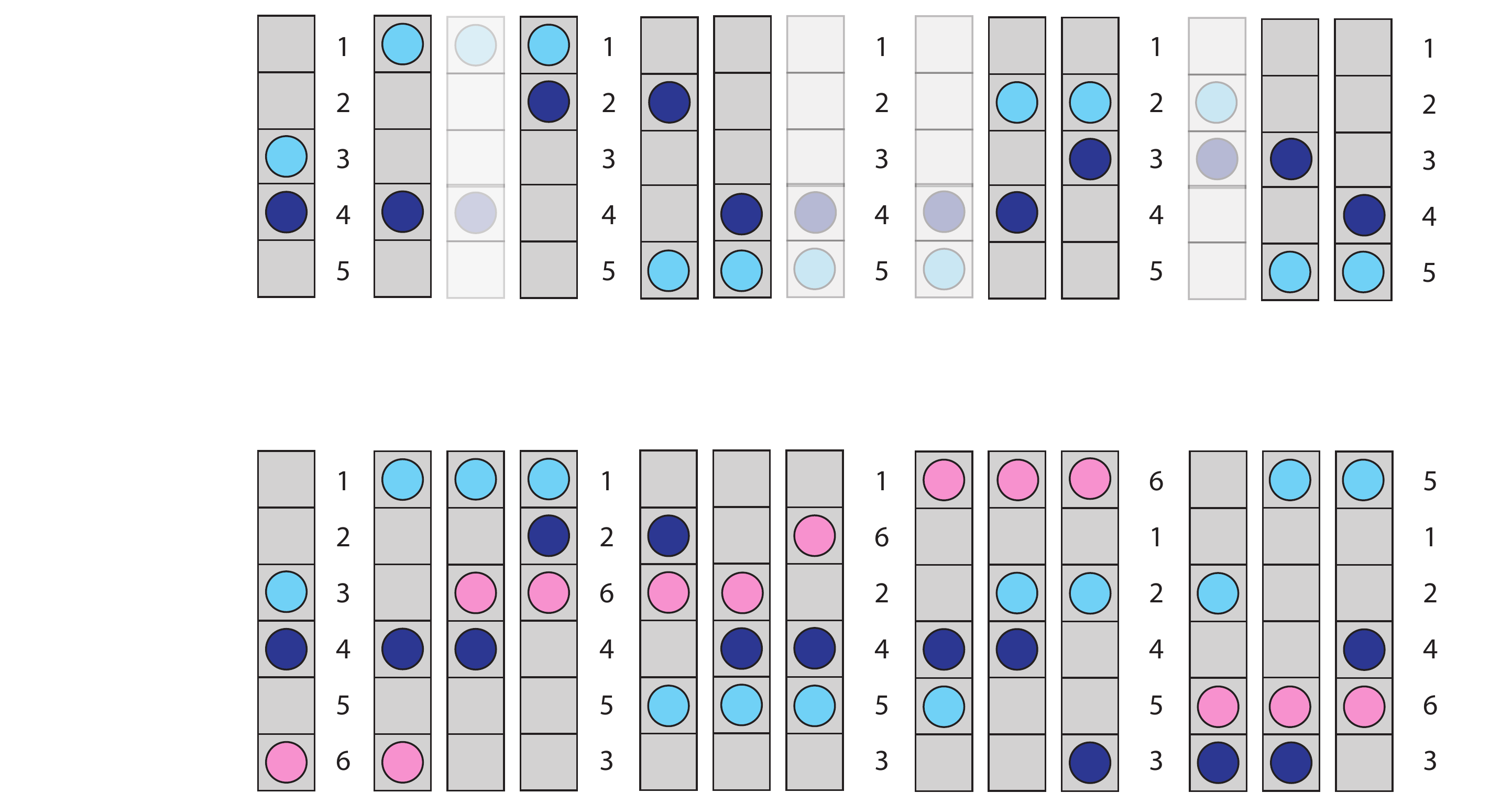}
  \rput(-7,7.8){\resizebox{1cm}{!}{$U_t$}}
  \rput(-7,2.6){\resizebox{1cm}{!}{$W_t$}}
  \rput(-5.9,5.3){\resizebox{0.8cm}{!}{$t=$}}
  \rput(-5.2,5.3){\resizebox{0.3cm}{!}{$0$}}
  \rput(-2.3,5.3){\resizebox{0.3cm}{!}{$1$}}
  \rput(0.7,5.3){\resizebox{0.3cm}{!}{$2$}}
  \rput(3.7,5.3){\resizebox{0.3cm}{!}{$3$}}
  \rput(6.7,5.3){\resizebox{0.3cm}{!}{$4$}}
  \caption{Above: $U_t$, the same SAC of $2$ walks on $K_5$ as in Figure~1. Faded are duplicates of previous steps used to synchronize with the extension below. Below: $W_t$, a POSAC of $3$ walks, under the partial order of the light blue walker walking before the dark blue one. The permutation is given at the end of every time unit while the order can be inferred from the diagram. Observe how the order of the blue walk change with respect to the extended pink walk between different time units. The rules is that the pink walk waits until
  his new place is clear and then movs. Also notice that the pink walk always ends his motion in place number $6$.}
\end{figure}
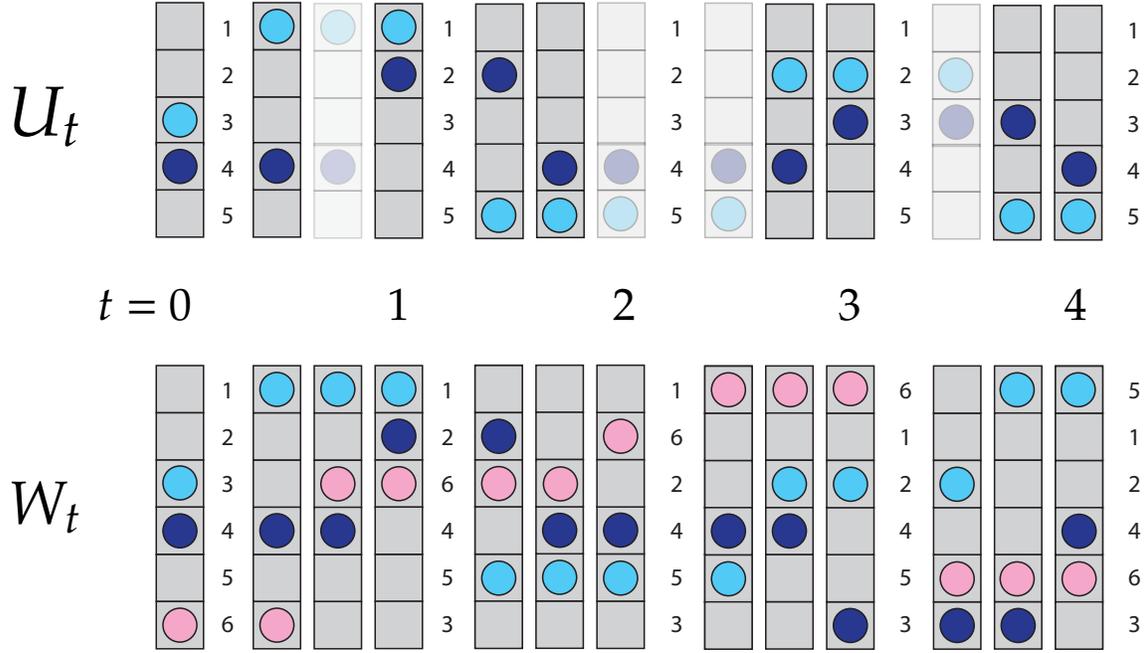

We begin by showing \eqref{eq: PSAC 1st}. Since the first $k$ walks of $\mathcal{W}^{N}_{k+1,R}$ are defined in
exactly the same way as these of $\mathcal{W}^{N}_{k}$, the proof that each of these walks performs a simple
random walk is identical to the proof of this fact for $\mathcal{W}^{N}_{k+1}$ and we omit it. The fact that
$\{W_t(k+1)\}_{t\in\Z}$ is a simple random walk is straightforward from fact that $W_t(k+1)=P_t(N)$ and from the definition of $P_t$.

Next let us show that $\mathcal{W}^{N}_{k+1,R}$ satisfies \eqref{eq: PSAC 2nd}. Observe that we have obtained the first part of
\eqref{eq: PSAC 2nd} in \ref{eq: first part uniqueness ext}.
For the second part, consider the event $$B_t^{i,j}=\{\sigma_t(i)<\sigma_t(j)\},$$ and write again $P'_t$ for the transposition $(N\ a_t)$.
For $i,j\in[k+1]$ we have
\begin{align*}
\Prob\Big(W_t(i)=W_{t-1}(j), B_t^{i,j}\Big)&=
\Prob\Big(P_{t}U_{t}(i)=P_{t-1}U_{t-1}(j),B_t^{i,j}\Big)=
\Prob\Big(P_{t-1}\circ P'_tU_{t}(i)=P_{t-1}U_{t-1}(j),B_t^{i,j}\Big)\\
&=\Prob\Big(P'_tU_{t}(i)=U_{t-1}(j),B_t^{i,j}\Big)=0,
\end{align*}\label{eq: SAC 2n pt 2}
following similar arguments to those used in \eqref{eq: SAC 2n pt 2}.

We thus are left with the case $k+1\in \{i,j\}$. However, if $i=k+1$ and $W_t(k+1)=W_{t-1}(j)$, then by
the definition of $\sigma_t$ we would have
$\sigma_t(j)=\sigma_t(k+1)=s_t(i)+1$ and $\sigma_t(i)=s_t(i)$. Thus
$$\Prob\Big(W_t(k+1)=W_{t-1}(j), B_t^{i,j}\Big)=0.$$

Finally consider the case $j=k+1$. If $B_t^{i,k+1}$ holds then, by the definition of $\sigma_t$ there must exist some $b\in[k]$ which satisfies
$B_t^{i,b}$ such that $W_{t-1}(b)=W_t(j)=P_t(N)$.
This $b$ satisfies $W_{t-1}(b)=P_{t-1}U_{t-1}(b)=P_t(N)$ and hence, by the definition of $P_t$, we have
$P_{t}U_{t-1}(b)=P_{t-1}(N)$.

We get that
\begin{align*}
\Prob\Big(W_t(i)=W_{t-1}(k+1), B_t^{i,j}\Big)&=\Prob\Big(W_t(i)=P_{t-1}(N), B_t^{i,j}\Big)
=\Prob\Big(\exists b\in[k]\ :\ P_tU_t(i)=P_{t}U_{t-1}(b), B_t^{i,b}\Big)\\&=
\Prob\Big(\exists b\in[k]\ :\  U_t(i)=U_{t-1}(b), B_t^{i,b}\Big)=
0.
\end{align*}
Where the last equality follows from the fact that $\mathcal{U}^{N-1}_{k,R}$ satisfies \eqref{eq: PSAC 2nd}.
Theorem~\ref{thm: main2} follows.

\section{Acknowledgments}

The author would like to thank David Wilson for introducing him to the problem and for useful discussions and to Omer Angel for his corrections and comments.

\end{document}